\documentclass[12pt]{amsart}

\usepackage{latexsym,amsfonts,amsmath,graphics,amssymb,multirow}

\newtheorem{remark}{Remark}
\newtheorem{example}{Example}

\newcommand{\dom}{{D}}

\newcommand{\bsDelta}{{\boldsymbol{\Delta}}}

\newcommand{\bsa}{{\boldsymbol{a}}}

\newcommand{\bsx}{{\boldsymbol{x}}}

\newcommand{\bsy}{{\boldsymbol{y}}}
\newcommand{\bsz}{{\boldsymbol{z}}}

\newcommand{\rd}{\mathrm{d}}

\newcommand{\bbE}{\mathbb{E}}

\newcommand{\bbR}{\mathbb{R}}
\newcommand{\bbZ}{\mathbb{Z}}
\newcommand{\bbN}{\mathbb{N}}
\newcommand{\calO}{\mathcal{O}}

\newcommand{\tr}{{\tt T}}

\allowdisplaybreaks

\usepackage[usenames]{color}
\definecolor{darkblue}{RGB}{40,0,200}
\definecolor{darkred}{RGB}{139,0,0}
\definecolor{darkgreen}{RGB}{0,100,0}
\definecolor{darkmagenta}{RGB}{139,0,139}

\begin{document}

\title{\scshape Fast QMC matrix-vector multiplication}

\author{Josef Dick, Frances Y. Kuo, Quoc T. Le Gia, Christoph Schwab}

\date{}
\maketitle

\begin{abstract}
Quasi-Monte Carlo (QMC) rules $1/N \sum_{n=0}^{N-1} f(\boldsymbol{y}_n A)$
can be used to approximate integrals of the form $\int_{[0,1]^s}
f(\boldsymbol{y} A) \,\mathrm{d} \boldsymbol{y}$, where $A$ is a matrix
and $\bsy$ is row vector. This type of integral arises for
example from the simulation of a normal distribution with a general
covariance matrix, from the approximation of the expectation value of
solutions of PDEs with random coefficients, or from applications from
statistics. In this paper we design QMC quadrature points
$\boldsymbol{y}_0, \ldots, \boldsymbol{y}_{N-1} \in [0,1]^s$ such that for
the matrix $Y = (\boldsymbol{y}_{0}^\top, \ldots,
\boldsymbol{y}_{N-1}^\top)^\top$ whose rows are the quadrature points, one
can use the fast Fourier transform to compute the matrix-vector product $Y
\boldsymbol{a}^\top$, $\boldsymbol{a} \in \mathbb{R}^s$, in $\mathcal{O}(N
\log N)$ operations and at most $s-1$ extra additions. The proposed method
can be applied to lattice rules, polynomial lattice rules and a certain
type of Korobov $p$-set.

The approach is illustrated computationally by three numerical
experiments. The first test considers the generation of points with normal
distribution and general covariance matrix, the second test applies QMC to
high-dimensional, affine-parametric, elliptic partial differential
equations with uniformly distributed random coefficients, and the third
test addresses Finite-Element discretizations of elliptic partial
differential equations with high-dimensional, log-normal random input
data. All numerical tests show a significant speed-up of the computation
times of the fast QMC matrix method compared to a conventional
implementation as the dimension becomes large.
\end{abstract}


{\bf MSC Class:} 65C05

\section{Introduction}

We are interested in numerical approximations of integrals of the form
\begin{equation} \label{eq:int}
\int_{U} f( \bsy A) \,\mu(\mathrm{d} \bsy),
\end{equation}
where the parameter domain $U$ is a subset of $\mathbb{R}^s$, $\mu$ is a
probability measure on $U$, $\boldsymbol{y}$ is a $1\times s$ row vector,
and $A$ is an $s \times t$ real matrix. Often we have $t = s$, but there
are also instances where $t$ is much larger than~$s$, see
Section~\ref{sec:app} below. We approximate these integrals by
equal-weight quadrature rules
\begin{equation}\label{eq_qmc}
  \frac{1}{N} \sum_{n=0}^{N-1} f(\boldsymbol{y}_n A),
\end{equation}
where $\bsy_0, \ldots, \bsy_{N-1} \in U$ are quadrature points which are
expressed again as row vectors (using row vectors merely simplifies our notation later on, it is not a necessity). We are interested in cases where the
computation of $ \boldsymbol{y}_n A$ for $n = 0, \ldots, N-1$ is a
significant factor in the computation of \eqref{eq_qmc} and where $N$ is
significantly smaller than $2^s$ (say $N \approx s^\kappa$ for some
$\kappa > 0$). The condition $N \ll 2^s$ is often naturally satisfied: for
instance if the dimension $s$ is very large, say $s > 100$, then the
number of points $N$ which can be used on current computers is much
smaller than $2^{100} \approx 10^{30}$; another instance arises for
example if the dimension is derived from a discretization or
approximation scheme where one needs to increase the dimension $s$
together with $N$ in order to reduce the discretization or approximation
error in a way such that $N \ll 2^s$. Examples where such situations arise
naturally are given in Section~\ref{sec:app}.

Returning to the approximation of \eqref{eq:int} by \eqref{eq_qmc}, one
concrete example of our setting arises from taking the expectation of some
quantity of interest with respect to the multivariate normal density with
a general covariance matrix $\Sigma\in\bbR^{s\times s}$,
\[
  \bbE[f] =
  \int_{\bbR^s} f(\bsz)\,\frac{\exp(-\frac{1}{2}\bsz\Sigma^{-1}\bsz^\tr)}{\sqrt{(2\pi)^s \det(\Sigma)}}\,\rd\bsz.
\]
Using a factorization $\Sigma = A^\tr A$ together with the substitution
$\bsz = \bsy A$, we arrive at the integral \eqref{eq:int}, with
$U = \bbR^s$, $s=t$, and with $\mu$ being the standard product Gaussian
measure. (If the mean for the multivariate normal density is nonzero then
a translation should be included in the substitution, but the general
principle remains the same.)
The method \eqref{eq_qmc} can be interpreted
as the simple \emph{Monte Carlo} approximation with $\bsy_0,
\ldots,\bsy_{N-1}\in\bbR^s$ being i.i.d.\ standard Gaussian random
vectors, and the computation of $\bsy_n A$ for $n=0,\ldots,N-1$ can be
interpreted as generating normally distributed points in $\mathbb{R}^s$
with the given covariance matrix $\Sigma$. The method \eqref{eq_qmc} can
also be interpreted as a \emph{quasi-Monte Carlo} (QMC) approximation with
$\bsy_n = \Phi^{-1}(\bsx_n)$ for $n=0,\ldots,N-1$, where
$\bsx_0,\cdots,\bsx_{N-1}\in [0,1]^s$ are deterministic QMC sample points,
and where $\Phi^{-1}: [0,1]\to\bbR$ denotes the inverse of the standard
normal distribution function and is applied component-wise to a vector. We
will consider this example in Subsection~\ref{num_exp:1}.

Several papers studied how one can obtain matrices $A$ in special
circumstances which allow a fast matrix-vector multiplication. In the
context of generating Brownian paths in mathematical finance, it is well
known that the ``standard construction'' (corresponding to the Cholesky
factorization of $\Sigma$) and the ``Brownian Bridge construction'' can be
done in $\calO(s)$ operations without explicitly carrying out the
matrix-vector multiplication (see e.g., \cite{Gla}), while the ``principal
components construction'' (corresponding to the eigenvalue decomposition
of $\Sigma$) can sometimes be carried out using Discrete Sine Transform in
$\calO(s\,\log s)$ operations (see e.g., \cite{GKSW08,Sche07}); other fast
orthogonal transform strategies can also be used (see e.g., \cite{Leo12}).
In the context of PDEs with random coefficients, it is known that
circulant embedding techniques can be applied in the generation of
stationary Gaussian random fields so that Fast Fourier Transform (FFT) can
be used (see e.g., \cite{DN97,GKNSS11}).

The approach in this paper differs from all of the above in that \emph{we
do not try to modify the matrix $A$, but rather modify the quadrature
points $\boldsymbol{y}_0, \ldots, \boldsymbol{y}_{N-1}$ to reduce the cost
of computing the $N$ matrix-vector products $ \boldsymbol{y}_n A$ for $n =
0, \ldots, N-1$. This implies that we do not require any structure in the
matrix $A$, and so our approach is applicable in general circumstances.}
For example, in some finance problems the payoff depends not only on the
Brownian paths but also on a basket of assets; our approach can be used to
speed up the matrix-vector multiplications with a factorization of the
covariance matrix among the assets. Another example arises from the
maximum likelihood estimation of generalised response models in
statistics: the change of variables strategy proposed in \cite{KDSWW08}
requires one to numerically compute the stationary point of the exponent
in the likelihood integrand function and the corresponding quadratic term
in the multivariate Taylor expansion; our approach can be used to speed up
the matrix-vector multiplications with a factorization of this numerically
computed Hessian matrix. Yet another important example arises from
parametric PDEs on high-dimensional parameter spaces, which appear in
computational uncertainty quantification; the presently proposed approach
can be used for both the so-called ``uniform'' and ``log-normal''
inputs (see e.g., \cite{GKNSSS13,GKNSS11,KSS12}). We will consider these PDE applications in
Subsections~\ref{num_exp:2} and~\ref{num_exp:3}.

To explain the idea behind our approach, we introduce the matrix
\begin{equation*}
Y = \begin{pmatrix} \boldsymbol{y}_0 \\ \vdots \\  \boldsymbol{y}_{N-1} \end{pmatrix} \in \mathbb{R}^{N \times s},
\end{equation*}
and we want to have a fast method to compute
\begin{equation*}
Y A = B = \begin{pmatrix} \boldsymbol{b}_0 \\ \vdots \\ \boldsymbol{b}_{N-1} \end{pmatrix}.
\end{equation*}
To compute \eqref{eq_qmc}, we propose to first compute the product $B = Y A$, store the matrix $B$, and then evaluate
\begin{equation*}
\frac{1}{N} \sum_{n=0}^{N-1} f(\boldsymbol{b}_n).
\end{equation*}
Thus this method requires $\mathcal{O}(N t)$ storage. In general, the
computation of $YA$ requires $\calO(Ns\,t)$ operations, and the quadrature
sum requires $\calO(N)$ operations. In the following we construct
quadrature points $\boldsymbol{y}_0, \ldots, \boldsymbol{y}_{N-1}$ for
which the matrix $Y$ permits a matrix-vector multiplication $Y\bsa$ in
$\mathcal{O}(N \log N)$ operations, where $\bsa$ can be any \emph{column}
of the matrix $A$. The computation of $Y A$ then reduces to
$\mathcal{O}(t\,N \log N)$ operations, instead of $\mathcal{O}(N s\,t)$
operations for the straight forward implementation. This leads to
significant speedup provided that $N$ is much smaller than $2^s$.

The basic idea of the proposed approach is to find quadrature point sets
$\{\boldsymbol{y}_0, \ldots, \boldsymbol{y}_{N-1}\}\in \mathbb{R}^s$ with
a specific ordering such that the matrix
$$Y' = \begin{pmatrix} \boldsymbol{y}_1 \\ \vdots \\ \boldsymbol{y}_{N-1}
\end{pmatrix} \in \mathbb{R}^{(N-1) \times s}$$ has a factorization of the
form
\begin{equation*}
Y' = Z P,
\end{equation*}
where $Z \in \mathbb{R}^{(N-1) \times (N-1)}$ is a circulant matrix and $P
\in \{0,1\}^{(N-1) \times s}$ is a matrix in which each column has at most
one value which is $1$ and with the remaining entries being $0$. The
special structure means that, for a given column vector $\bsa$, the column vector $\bsa' = P \boldsymbol{a}$
can be obtained in at most $\mathcal{O}(N)$ operations, and the
matrix-vector multiplication $Z \bsa'$ can be computed in $\mathcal{O}(N
\log N)$ operations using FFT. On the other hand, the computation of
$\boldsymbol{y}_0 \boldsymbol{a}$ requires at most $s-1$ additions and $1$
multiplication. The vector $\bsy_0$ is separated out because typical QMC
methods would lead to $\bsy_0$ being a constant vector. If
$\boldsymbol{y}_0 = (0, \ldots, 0)$, then no extra computation is necessary.

In Section~\ref{sec_fast} we consider two important classes of QMC
point sets whose structure facilitates the use of the presently proposed
acceleration:
\begin{itemize}
\item point sets derived from lattice rules,
\item the union of all Korobov lattice point sets (which is one class of ``Korobov $p$-sets'').
\end{itemize}

The same strategy can be applied also to polynomial lattice rules
where the modulus is a primitive polynomial over a finite field $\mathbb{F}_b$ of order $b$, and to the union of all
Korobov polynomial lattice rules.

Note that lattice rules and polynomial lattice rules can yield a
convergence rate close to $\calO(N^{-1})$ for sufficiently smooth
integrands, with the implied constant independent of the
integration-dimension $s$ under appropriate conditions on the integrand
function and the underlying function space setting, see e.g.,
\cite{DKS13}. The union of Korobov lattice point sets on the other
hand achieves a convergence rate of $\calO(N^{-1/2})$ for a much larger
class of functions, and dimension independent error bounds can be obtained
with significantly weaker assumptions \cite{D14, DP14}. Thus, when the
integrand is not smooth enough for lattice rules and polynomial lattice
rules, the union of Korobov lattice point sets 
can be a good substitute for the simple Monte Carlo method so that the
fast computation approach of this paper can be exploited.

To illustrate our method and to investigate numerically for which
parameter ranges the improvements in the computational cost are visible, we consider
three applications in Section~\ref{sec:app}. In
Subsection~\ref{num_exp:1} we generate normally distributed points with a
general covariance matrix. In Subsections~\ref{num_exp:2}
and~\ref{num_exp:3} we consider PDEs with random coefficients in the
uniform case and log-normal case, respectively. The numerical results in Section~\ref{sec:num} show that our method
is significantly faster whenever the dimension becomes large.

\section{Fast QMC matrix-vector multiplication}\label{sec_fast}

We explain the fast method for lattice point sets and the union of
all Korobov lattice point sets. The basic idea also applies to polynomial
lattice point sets and the union of all Korobov polynomial lattice point
sets.

\subsection{Fast matrix-vector multiplication for lattice point sets}\label{sec_fast_lattice}

Our approach is very similar to the method used in \cite{NC06a} for the
fast component-by-component construction of the generating vector for
lattice rules, however, we apply it now to the matrix vector multiplication $Y \bsa$ rather then the component-by-component construction. For simplicity, we confine the exposition to cases where
the number of points is a prime. Based on the presently developed ideas, the general case can be handled
analogously with the method from \cite{NC06b}.

Let $N$ be a prime number, let $\mathbb{Z}_N = \{0, 1, \ldots, N-1\}$, and
let $\mathbb{Z}_N^\ast = \{1,2, \ldots, N-1\}$.

A lattice point set with generator $(g_1, g_2, \ldots, g_s) \in
(\mathbb{Z}^\ast_N)^s$ is of the form
\begin{equation*}
\left(\left\{\frac{n g_1}{N}\right\}, \left\{ \frac{n g_2}{N} \right\},
\ldots, \left\{\frac{n g_{s}}{N}\right\} \right) \quad \mbox{for } n = 0, 1, \ldots, N-1,
\end{equation*}
where for nonnegative real numbers $x$ we denote by $\{x\} = x - \lfloor x
\rfloor$ the fractional part of $x$.

Let $\beta$ be a primitive element of the multiplicative group
$\mathbb{Z}_N^\ast$, i.e., we have $\{\beta^{\,k} \bmod N :
k=1,2,\ldots,N-1\} = \bbZ_N^*$. As is well-known $\beta^{N-1} \equiv \beta^0
\equiv 1 \bmod N$. Moreover, its multiplicative inverse $\beta^{-1}
\in \mathbb{Z}_N^*$ is also a primitive element. We write each component
of the generating vector $(g_1,g_2\ldots,g_s)$ as
\[
  g_j \equiv \beta^{\,c_j-1} \bmod N \quad\mbox{for some}\quad 1 \le c_j \le N-1.
\]
Note that the fast component-by-component algorithm of \cite{NC06a} for
constructing the generating vector computes the values $c_j$ as a
by-product, and hence no additional computation is needed to obtain the
values $c_j$ in this case.

Clearly, the ordering of the QMC points does not affect the quadrature sum.
We now specify a particular (unconventional) ordering which allows fast
matrix-vector multiplications. We define $\boldsymbol{x}_0 = (0, \ldots,
0)$, and for $n=1,2,\ldots,N-1$ we define
\begin{align*}
 \boldsymbol{x}_n
 &= \left(\left\{\frac{\beta^{-(n-1)} g_1 }{N}\right\},
 \left\{ \frac{\beta^{-(n-1)} g_2}{N} \right\}, \ldots,
 \left\{\frac{\beta^{-(n-1)} g_s}{N}\right\} \right) \\
 &= \left(\left\{\frac{\beta^{-(n-1)} \beta^{\,c_1-1} }{N}\right\},
 \left\{ \frac{\beta^{-(n-1)} \beta^{\,c_2-1}}{N} \right\}, \ldots,
 \left\{\frac{ \beta^{-(n-1)} \beta^{\,c_s-1}}{N}\right\} \right) \\
 &=  \left(\left\{\frac{ \beta^{\,c_1-n} }{N}\right\},
 \left\{ \frac{ \beta^{\,c_2-n} }{N} \right\}, \ldots, \left\{\frac{ \beta^{\,c_s-n} }{N}\right\} \right).
\end{align*}
In essence, we have changed the ordering by substituting the conventional
index $n$ with $ \beta^{-(n-1)}$ and replacing each generating vector
component $g_j$ by $ \beta^{\,c_j-1}$.

The quadrature points we consider in \eqref{eq_qmc} are now given by
\begin{align*}
 \boldsymbol{y}_n
 &= \varphi(\boldsymbol{x}_n) \\
 &= (\varphi(x_{n,1}), \varphi(x_{n,2}), \ldots,\varphi(x_{n,s}))
 \quad \mbox{for } n = 0, 1, \ldots, N-1,
\end{align*}
where \emph{we apply the same univariate transformation $\varphi:[0,1] \to
\mathbb{R}$ to every component of every point.} One example for such a
transformation is $\varphi(x) =\Phi^{-1}(x)$, the inverse of the
cumulative normal distribution function; this maps the points from
$[0,1]^s$ to $\bbR^s$ as we already discussed in the introduction. Another
example is $\varphi(x) = 1 - |2x-1|$, the tent transform; results for lattice rules usually apply to periodic functions, applying the tent transform yields similar results for non-periodic functions, see
\cite{DNP13}. The case where $\varphi(x) = x$ is included as a special
case.

We discuss now the multiplication of the matrix $Y$ with a column vector
$\boldsymbol{a} \in \mathbb{R}^s$. Since $\boldsymbol{y}_0 = (\varphi(0),
\varphi(0), \ldots, \varphi(0) )$ we have $$\boldsymbol{y}_0
\boldsymbol{a}  = \varphi(0) \sum_{j=1}^s a_j.$$ In particular, if
$\varphi$ is the identity mapping then $\boldsymbol{y}_0 \boldsymbol{a} =
0$.
Thus the first component can be computed using at most $s-1$ additions
and $1$ multiplication. We consider now the remaining matrix
\begin{equation*}
Y' = \begin{pmatrix} \boldsymbol{y}_1 \\ \vdots \\ \boldsymbol{y}_{N-1} \end{pmatrix}.
\end{equation*}
In the following we show that $Y'$ can be written as a product of a
circulant matrix $Z$ and a matrix $P$ in which $N-1$ entries are $1$ and
the remaining entries are $0$.

Recall that $ \beta$ is a primitive element of $\mathbb{Z}_N^\ast$.
For $k \in \mathbb{Z}$ let
\begin{equation*}
z_{k} = \varphi\left( \left\{\frac{ \beta^{k}}{N} \right\} \right).
\end{equation*}
Then we have $z_{k} = z_{k + \ell (N-1)}$ for all $ \ell \in \mathbb{Z}$.
Let
\begin{equation*}
 Z  = \begin{pmatrix}
 z_0 & z_1 & z_2 & \ldots & z_{N-3} & z_{N-2} \\
 z_{N-2} & z_0 & z_1 & \ddots & \ddots & z_{N-3} \\
 z_{N-3} & z_{N-2} & z_0 & \ddots & \ddots & \vdots \\
 \vdots & \ddots & \ddots & \ddots & \ddots & \vdots \\
 z_2 & \ddots &  \ddots & \ddots & z_0 & z_1 \\
 z_1 & z_2 & \ldots & \ldots & z_{N-2} & z_0
 \end{pmatrix}.
\end{equation*}
We define the matrix $P = (p_{k,j})_{1 \le k \le N-1, 1 \le j \le s} \in
\{0, 1\}^{(N-1) \times s}$ by
\begin{equation*}
  p_{k,j}
  = \begin{cases}
  1 & \mbox{if } k = c_j, \\
  0 & \mbox{otherwise}.
  \end{cases}
\end{equation*}
Each column of the matrix $P$ contains exactly one element $1$, with the
remaining elements being $0$. It is now elementary to check that
\begin{equation}\label{ZPa}
Y' = Z P.
\end{equation}
Note that the matrix $Z$ is exactly the same as the matrix used in the
fast component-by-component algorithm of \cite{NC06a}. In effect, the
matrix $P$ specifies which columns of $Z$ to select (namely, the $c_1$-th,
the $c_2$-th, \ldots, and the $c_s$-th) to recover $Y'$.

Let $\boldsymbol{a} \in \mathbb{R}^{s}$ be any column vector. Then $\bsa'
= P \boldsymbol{a}$ can be obtained in at most $\mathcal{O}(N)$ operations
due to the special structure of $P$, and the matrix-vector multiplication
$Z \bsa'$ can be computed in $\mathcal{O}(N \log N)$ operations using FFT (see \cite{FFTW05})
since $Z$ is circulant. Thus $Y' \boldsymbol{a}$ can be computed in
$\mathcal{O}(N \log N)$ operations, and hence the matrix-vector
multiplication $Y \boldsymbol{a}$ can be carried out using $\mathcal{O}(N
\log N)$ operations plus at most $s-1$ additions.

We remark that the formula \eqref{ZPa} can also be used to generate the
matrix $Y'$, i.e., to generate the quadrature points in a fast way.
Further, if one wants to store the point set, i.e., matrix $Y$, one can
simply store the primitive root $\beta$ and the $s$ numbers $c_1,\ldots,c_s$.

The case where $N$ is not a prime number can be treated as in
\cite{NC06b}.

We finish this subsection with a simple example to illustrate the idea.

\begin{example}
Let $s=3$, $N=7$, and $(g_1,g_2,g_3) = (1,5,3)$. A primitive root for
$\bbZ_7^*$ is $ \beta=3$, with multiplicative inverse $\beta^{-1}=5$. We
have
\begin{align*}
  g_1 &= 1 = 3^{1-1} \bmod 7 \quad\implies\quad c_1 = 1, \\
  g_2 &= 5 = 3^{6-1} \bmod 7 \quad\implies\quad c_2 = 6, \\
  g_3 &= 3 = 3^{2-1} \bmod 7 \quad\implies\quad c_3 = 2.
\end{align*}
The conventional ordering of the points and the new ordering are
\begin{align*}
  \begin{cases}
  (0, 0,0), \\
  (\tfrac{1}{7}, \tfrac{5}{7}, \tfrac{3}{7}), \\
  (\tfrac{2}{7}, \tfrac{3}{7}, \tfrac{6}{7}), \\
  (\tfrac{3}{7}, \tfrac{1}{7}, \tfrac{2}{7}), \\
  (\tfrac{4}{7}, \tfrac{6}{7}, \tfrac{5}{7}), \\
  (\tfrac{5}{7}, \tfrac{4}{7}, \tfrac{1}{7}), \\
  (\tfrac{6}{7}, \tfrac{2}{7}, \tfrac{4}{7}),
  \end{cases}
  \;\mbox{versus}\quad
  \begin{cases}
  \bsx_0 = (0, 0,0), \\
  \bsx_1 =
  (\{\frac{3^{1-1}}{7}\},\{\frac{3^{6-1}}{7} \},\{\frac{3^{2-1}}{7}\} )
  = (\frac{1}{7}, \frac{5}{7}, \frac{3}{7}), \\
  \bsx_2 =
  (\{\frac{3^{1-2}}{7}\},\{\frac{3^{6-2}}{7} \},\{\frac{3^{2-2}}{7}\} )
  = (\frac{5}{7}, \frac{4}{7}, \frac{1}{7}), \\
  \bsx_3 =
  (\{\frac{3^{1-3}}{7}\},\{\frac{3^{6-3}}{7} \},\{\frac{3^{2-3}}{7}\} )
  = (\frac{4}{7}, \frac{6}{7}, \frac{5}{7}), \\
  \bsx_4 =
  (\{\frac{3^{1-4}}{7}\},\{\frac{3^{6-4}}{7} \},\{\frac{3^{2-4}}{7}\} )
  = (\frac{6}{7}, \frac{2}{7}, \frac{4}{7}), \\
  \bsx_5 =
  (\{\frac{3^{1-5}}{7}\},\{\frac{3^{6-5}}{7} \},\{\frac{3^{2-5}}{7}\} )
  = (\frac{2}{7}, \frac{3}{7}, \frac{6}{7}), \\
  \bsx_6 =
  (\{\frac{3^{1-6}}{7}\},\{\frac{3^{6-6}}{7} \},\{\frac{3^{2-6}}{7}\} )
  = (\frac{3}{7}, \frac{1}{7}, \frac{2}{7}).
  \end{cases}
\end{align*}
It is easy to see that indeed
\begin{align*}
  \underset{Y'}{
  \begin{pmatrix}
  \varphi(\bsx_1) \vphantom{\frac{1}{7})}\vspace{0.1cm} \\
  \varphi(\bsx_2) \vphantom{\frac{1}{7})}\vspace{0.1cm} \\
  \varphi(\bsx_3) \vphantom{\frac{1}{7})}\vspace{0.1cm} \\
  \varphi(\bsx_4) \vphantom{\frac{1}{7})}\vspace{0.1cm} \\
  \varphi(\bsx_5) \vphantom{\frac{1}{7})}\vspace{0.1cm} \\
  \varphi(\bsx_6)
  \end{pmatrix}
  }
  =
  \underset{Z}{
  \begin{pmatrix}
  \varphi(\frac{1}{7}) & \varphi(\frac{3}{7}) & \varphi(\frac{2}{7}) & \varphi(\frac{6}{7}) & \varphi(\frac{4}{7}) & \varphi(\frac{5}{7}) \vspace{0.1cm} \\
  \varphi(\frac{5}{7}) & \varphi(\frac{1}{7}) & \varphi(\frac{3}{7}) & \varphi(\frac{2}{7}) & \varphi(\frac{6}{7}) & \varphi(\frac{4}{7}) \vspace{0.1cm} \\
  \varphi(\frac{4}{7}) & \varphi(\frac{5}{7}) & \varphi(\frac{1}{7}) & \varphi(\frac{3}{7}) & \varphi(\frac{2}{7}) & \varphi(\frac{6}{7}) \vspace{0.1cm} \\
  \varphi(\frac{6}{7}) & \varphi(\frac{4}{7}) & \varphi(\frac{5}{7}) & \varphi(\frac{1}{7}) & \varphi(\frac{3}{7}) & \varphi(\frac{2}{7}) \vspace{0.1cm} \\
  \varphi(\frac{2}{7}) & \varphi(\frac{6}{7}) & \varphi(\frac{4}{7}) & \varphi(\frac{5}{7}) & \varphi(\frac{1}{7}) & \varphi(\frac{3}{7}) \vspace{0.1cm} \\
  \varphi(\frac{3}{7}) & \varphi(\frac{2}{7}) & \varphi(\frac{6}{7}) & \varphi(\frac{4}{7}) & \varphi(\frac{5}{7}) & \varphi(\frac{1}{7})
  \end{pmatrix}
  }
  \underset{P}{
  \begin{pmatrix}
  1 & 0 & 0 \vphantom{\frac{1}{7})}\vspace{0.1cm} \\
  0 & 0 & 1 \vphantom{\frac{1}{7})}\vspace{0.1cm} \\
  0 & 0 & 0 \vphantom{\frac{1}{7})}\vspace{0.1cm} \\
  0 & 0 & 0 \vphantom{\frac{1}{7})}\vspace{0.1cm} \\
  0 & 0 & 0 \vphantom{\frac{1}{7})}\vspace{0.1cm} \\
  0 & 1 & 0
  \end{pmatrix}
  }.
\end{align*}
The matrix $P$ specifies that we select the first, the sixth, and the
second columns of $Z$, as indicated by the values of $c_1, c_2, c_3$, to
recover $Y'$.
\end{example}

\begin{remark}
The method discussed above works in the same way for polynomial lattice
rules over a finite field $\mathbb{F}_b$ of order $b$.
\end{remark}

\begin{remark}
The method does \emph{not} work when we apply general randomization
techniques such as ``shifting'' for lattice rules or ``scrambling'' for
polynomial lattice rules. This is because the corresponding transformation
$\varphi$ in the mapping $\bsy_n = \varphi(\bsx_n)$ fails to be the same
mapping in all coordinate directions. If we were to restrict all random
shifts to be of the form $\bsDelta = (\Delta,\ldots,\Delta) \in [0,1]^s$
in the case of lattice rules then the method would work.
\end{remark}

\begin{remark}
Higher order polynomial lattice rules, which have been introduced in \cite{DP07}, also fit into the structure used in this subsection since they can be viewed as the first $b^m$ points of a polynomial lattice point sets with $b^{m \alpha}$ points, where $\alpha \in \mathbb{N}$ denotes the smoothness. Here $\alpha = 1$ corresponds to the classical polynomial lattice rules. However, if we use the method from this paper, then the matrix vector multiplication uses the full $b^{m \alpha}$ points, which means the matrix vector multiplication requires $\mathcal{O}(b^{m\alpha} m \alpha)$ operations, instead of $\mathcal{O}(b^m s)$ operations for a straightforward implementation. Thus this method is only advantageous if $b^{m \alpha} m \alpha \ll b^m s$. For $\alpha \ge 2$ this implies that $b^m \ll s$, which usually does not hold.
\end{remark}

\subsection{The union of all Korobov lattice point sets}

Hua and Wang~\cite[Section~4.3]{HW81} studied the point set
\begin{equation*}
 \left(\left\{ \frac{n g^0}{K} \right\}, \left\{ \frac{n g^1}{K} \right\},
 \ldots, \left\{\frac{n g^{s-1}}{K} \right\} \right), \quad \mbox{for } n, g = 1, 2 \ldots, K-1,
\end{equation*}
where $K$ is a prime number. The number of points is $N = (K-1)^2$. This
is essentially the union of all Korobov lattice point sets. (Note that Hua
and Wang also included the cases $n = 0$ or $g=0$, or both $n=g=0$, but
these only yield the zero vector.)

This point set achieves only a rate of convergence of the weighted
star-discrepancy of $\mathcal{O}(N^{-1/2 + \delta})$ for any $\delta > 0$,
however, the dependence on the dimension of the weighted star-discrepancy
is better than what is known for lattice point sets or polynomial lattice
point sets in some circumstances, see \cite{DP14} for more details.

We now specify a particular ordering of the points to allow fast
matrix-vector multiplications. Let $\beta$ be a primitive element in
$\mathbb{Z}_K^*$. As in Subsection~\ref{sec_fast_lattice}, we replace the
index $n$ in the conventional ordering by $\beta^{-(n-1)}$, and
similarly we replace the index $g$ by $\beta^{(g-1)}$. That is, for
$n, g = 1, 2, \ldots, K-1$, we define
\begin{align*}
\boldsymbol{x}_{n,g} &= \left(\left\{ \frac{\beta^{-(n-1)} \beta^{0(g-1)}}{K} \right\},
\left\{ \frac{ \beta^{-(n-1)} \beta^{1(g-1)}}{K} \right\}, \left\{\frac{ \beta^{-(n-1)} \beta^{2(g-1)}}{K} \right\},
\ldots,  \right.\\
&\qquad\qquad\qquad\qquad\qquad\qquad\qquad\qquad\qquad \left.
 \left\{\frac{ \beta^{-(n-1)} \beta^{(s-1)(g-1)}}{K} \right\} \right) \\
 &= \left(\left\{ \frac{ \beta^{c_{g,1}-n}}{K} \right\},
\left\{ \frac{ \beta^{c_{g,2}-n}}{K} \right\}, \left\{\frac{ \beta^{c_{g,3}-n}}{K} \right\},
\ldots, \left\{\frac{ \beta^{c_{g,s}-n}}{K} \right\} \right),
\end{align*}
with
\[
  c_{g,j} = (j-1)(g-1)+1 \quad\bmod{(K-1)} \qquad\mbox{for}\quad j=1,\ldots,s.
\]
We also define
\begin{equation*}
\boldsymbol{y}_{n,g} = \varphi(\boldsymbol{x}_{n,g}).
\end{equation*}

Finally we define the matrix
\begin{equation}\label{Y_pset}
 Y' = \begin{pmatrix} Y_1' \\ Y_2' \\ \vdots \\ Y_{K-1}' \end{pmatrix},
 \quad\mbox{with}\quad
 Y_g' = \begin{pmatrix} \boldsymbol{y}_{1, g} \\ \boldsymbol{y}_{2, g} \\ \vdots \\ \boldsymbol{y}_{K-1, g}  \end{pmatrix}
 \;\mbox{for}\; g=1,\ldots,K-1.
\end{equation}
For the matrices $Y_g'$ we can apply the method from
Subsection~\ref{sec_fast_lattice} to write it as $Y_g' = Z P_g$ using the
values of $c_{g,1},\ldots,c_{g,s}$ so that a matrix-vector multiplication
can be computed in at most $\mathcal{O}(K \log K)$ operations.
Thus one matrix-vector product for the matrix $Y'$ 
can be evaluated in $\mathcal{O}(N \log N)$ operations.

\begin{remark}
The same strategy can be applied to the union of all Korobov polynomial
lattice point sets.
\end{remark}

\section{Applications} \label{sec:app}

\subsection{Generation of normally distributed points with general covariance matrix}
\label{num_exp:1}

In many applications one requires realizations of random variables in
$\mathbb{R}^s$ which are normally distributed
$\mathcal{N}(\boldsymbol{\mu}, \Sigma)$ with mean $\boldsymbol{\mu} =
(\mu_1, \mu_2, \ldots, \mu_s)$ and covariance matrix $\Sigma \in
\mathbb{R}^{s \times s}$. An algorithm to generate such random variables
is described for instance in \cite[Section~11.1.6]{HLD} and works the
following way. Let $A \in \mathbb{R}^{s \times s}$ be such that $A^\top A
= \Sigma$; for example, $A$ can be the upper triangular matrix in the
Cholesky decomposition of $\Sigma$. To generate a point $(Z_1, Z_2,
\ldots, Z_s) \sim \mathcal{N}(\boldsymbol{\mu}, \Sigma)$, one generates
i.i.d.\ standard normal random variables $Y_1, Y_2, \ldots, Y_s \sim
\mathcal{N}(0,1)$ with mean $0$ and variance $1$ and then computes $
  (Z_1, Z_2, \ldots, Z_s) =
  (Y_1, Y_2,\ldots, Y_s) A + (\mu_1, \mu_2, \ldots, \mu_s).
$

As we already outlined in the introduction, this procedure can be
implemented in the following way using deterministic QMC point sets. Let
$\boldsymbol{x}_0, \boldsymbol{x}_1, \ldots, \boldsymbol{x}_{N-1} \in
[0,1]^s$ be a set of quadrature points as described in
Section~\ref{sec_fast}. Let $\Phi^{-1}$ be the inverse of the cumulative
normal distribution function. Set
\begin{equation*}
\boldsymbol{y}_n = \Phi^{-1}(\boldsymbol{x}_n)
\end{equation*}
and compute
\begin{equation}\label{Ay_prod}
\boldsymbol{z}_n = \boldsymbol{y}_n A + \boldsymbol{\mu}.
\end{equation}
Note that we do not assume any structure in the matrix $A$. This is
contrary to a number of scenarios in e.g., mathematical finance; see our
discussion in the introduction.

\subsection{Partial differential equations with ``uniform'' random coefficients}\label{num_exp:2}

The matrix-vector multiplication also arises in applications of QMC
for approximating linear functionals of solutions of PDEs with random
coefficients, see e.g., \cite{KSS12}.

A prototypical class of countably parametric, elliptic boundary value
problems reads as
\begin{align}
  -\nabla \cdot (\mathfrak{a}(\vec{x}, \boldsymbol{y}) \nabla u(\vec{x}, \boldsymbol{y}))
  &= g(\vec{x}), && \vec{x} \in \dom\subset \mathbb{R}^d,\;
     \boldsymbol{y} \in [-\tfrac{1}{2}, \tfrac{1}{2}]^{\mathbb{N}}, \label{PDE_uniform} \\
  u(\vec{x}, \boldsymbol{y}) &=  0, && \vec{x} \in \partial \dom, \label{PDE_boundary}
\end{align}
where $\vec{x}$ is a vector in the convex, physical domain $\dom\subseteq\bbR^d$,
and $\boldsymbol{y} = (y_1, y_2, \ldots)$ is a parametric sequence in
$[-\tfrac{1}{2}, \tfrac{1}{2}]^\bbN$, with $y_j$ being uniformly
distributed on $[-\tfrac{1}{2}, \tfrac{1}{2}]$. Here, we distinguish
notationally between a vector $\vec{x}$ in the spatial domain $\dom$ and a vector
$\bsy$ in the parametric domain $[-\tfrac{1}{2}, \tfrac{1}{2}]^{\mathbb{N}}$. The coefficient
$\mathfrak{a}(\vec{x},\bsy)$ depends on the parameter sequence $\bsy$ in an
affine manner, ie.
\begin{equation} \label{eq:axy}
\mathfrak{a}(\vec{x}, \boldsymbol{y})
=
\psi_0(\vec{x}) + \sum_{j=1}^\infty y_j\, \psi_j(\vec{x})
\;
\end{equation}
In \eqref{eq:axy}, the sequence of functions $\psi_j(\vec{x})$ for $j\ge
1$ is assumed to decay as $j\to\infty$ such that $ \sum_{j\geq 1} \|
\psi_j \|_{L^\infty(\dom)}^p < \infty$ for some $0 < p \leq 1$ and that $\sum_{j \ge 1} \|\nabla \psi_j\|_{L^\infty(D)} < \infty$. For more background, necessary assumptions and theoretical results see
\cite{KSS12}.

In \cite{KSS12} the aim is to approximate the expected value
$\mathbb{E}(G(u))$ with respect to the random sequence $\boldsymbol{y}$ of
a linear functional $G$ of the solution $u$ of the PDE. The algorithm
truncates the infinite sum in \eqref{eq:axy} after $s$ terms (i.e. set
$y_j = 0$ for $j > s$), solves the truncated problem using a (piecewise
linear) finite element method with $M$ mesh points, and then approximates
the $s$-dimensional integral using an $N$-point QMC rule. The resulting
three sources of errors, namely, the truncation error, the finite element
error, and the quadrature error, need to be balanced. For instance,
in the case that $\sum_{j=1}^\infty \|\psi_j\|_{L_\infty(\dom)}^{2/3} <
\infty$ and that $g,G\in L^2(\dom)$, \cite[Theorem~8.1]{KSS12} yields for
continuous, piecewise linear Finite Element discretizations of \eqref{PDE_uniform} -- \eqref{eq:axy}  in $\dom$ on
quasiuniform meshes that we should choose
\begin{equation} \label{PDE_Ns}
  N \asymp s \asymp M^{2/d},
\end{equation}
where $\asymp$ indicates that the terms should be of the same order. In
general we have $N \asymp s^\kappa$ for some small $\kappa > 0$ (i.e., $N
\ll 2^s$). Thus the fast method of this paper can be advantageous (see the numerical results in Section~\ref{sec:num}).

 Let $\phi_1, \phi_2, \ldots, \phi_M$ be
a basis for the finite element space $V_M$. For each $0\le n < N$, let
$u_M^{(n)} = \sum_{k=1}^M \widehat{u}_k^{(n)} \phi_k$ be the finite
element approximation of the solution $u$ of the PDE given the parameter
$\bsy_n\in [-\frac{1}{2},\frac{1}{2}]^s$. Then for each $0\le n < N$ we
solve the linear system
\begin{equation}\label{PDE_lin_sys}
  (\widehat{u}_1^{(n)}, \widehat{u}_2^{(n)}, \ldots,\widehat{u}_M^{(n)})\, B(\boldsymbol{y}_n)
  = (\widehat{g}_1, \widehat{g}_2, \ldots,\widehat{g}_M)
\end{equation}
for $(\widehat{u}_1^{(n)}, \widehat{u}_2^{(n)},
\ldots,\widehat{u}_M^{(n)})$, where $\widehat{g}_k = \int_\dom
g(\vec{x})\,\phi_k(\vec{x})\,\rd\vec{x}$ for $k=1,\ldots,M$, and where the
symmetric stiffness matrix $B(\bsy_n) = (b_{n,k,\ell})_{k,\ell}$ depends
on $\bsy_n$ and has entries
\begin{equation} \label{a_approx}
  b_{n,k,\ell} =
  \int_\dom \mathfrak{a}(\vec{x},\bsy_n)\,\nabla \phi_k(\vec{x})\cdot \nabla \phi_\ell(\vec{x})\,\rd\vec{x},
  \qquad 1\le k,\ell \le M.
\end{equation}
Note that, with a slight abuse of notation, $\mathfrak{a}(\vec{x},\bsy_n)$
is given by \eqref{eq:axy} but truncated to $s$ terms. The expected value
of the solution can then be approximated by
\begin{equation}\label{eq_sol_PDE}
\overline{u}(\vec{x}) = \frac{1}{N} \sum_{n=0}^{N-1} \sum_{k=1}^{M-1} \widehat{u}^{(n)}_k \phi_k(\vec{x}).
\end{equation}

Due to the linear structure in \eqref{eq:axy} for the uniform case, we can
write
\[
  b_{n,k,\ell} = a_{0,k,\ell} +  \sum_{j=1}^s y_{n,j}\, a_{j,k,\ell},
  \qquad 0\le n < N,\quad 1\le k,\ell \le M,
\]
where $y_{n,j}$ denotes the $j$th component of the $n$th point $\bsy_n$,
and
\begin{equation*}
  a_{j, k, \ell} = \int_D \psi_j(\vec{x}) \nabla \phi_k(\vec{x}) \cdot \nabla \phi_\ell(\vec{x}) \,\mathrm{d} \vec{x},
  \qquad 1\le k,\ell \le M, \quad j\ge 0.
\end{equation*}
In the standard approach, one defines the symmetric matrices $A_j =
(a_{j,k,\ell})_{1 \le k, \ell \le M}$ and sets
\begin{equation*}
B(\boldsymbol{y}_n ) = A_0 + \sum_{j=1}^s y_{n,j}\, A_j,
\qquad 0\le n < N.
\end{equation*}
Note that $A_j$ is usually sparse, with only $\calO(M)$ nonzero entries in the same position, depending only on the relative supports of the basis functions $\phi_k$, which are thus in particular independent of $j$.
The cost for computing $B(\bsy_n)$ for all $n=0,\ldots,N-1$ is therefore
$\calO(M N s)$ operations.

The fast QMC matrix-vector approach is implemented as follows: let $Y =
(\boldsymbol{y}_n)_{0 \le n < N}$ be the matrix whose rows are the
quadrature points of the QMC rule, and let $\boldsymbol{a}_{k, \ell} =
(a_{1, k, \ell}, \ldots, a_{s,k,\ell})^\top$ and $\boldsymbol{b}_{k, \ell
} = (b_{0,k,\ell},\ldots,b_{N-1,k,\ell})^\top$. Then compute
\begin{equation}\label{Ya}
 \boldsymbol{b}_{k, \ell }
 = (a_{0,k,\ell}, \ldots, a_{0, k, \ell})^\top + Y  \boldsymbol{a}_{k, \ell}
 \quad \mbox{for all } 1 \le k, \ell \le M,
\end{equation}
where each matrix-vector multiplication $Y \boldsymbol{a}_{k, \ell}$
should be done using the approach of the previous section. Since only
$\mathcal{O}(M)$ vectors $\bsa_{k, \ell}$ are nonzero, this approach for
obtaining $B(\bsy_n)$ for all $n=0,\ldots,N-1$ therefore requires only
$\calO(M\, N \log N)$ operations.

The improvement in the computational cost is that we replaced a factor of
$\calO(s)$ by $\calO(\log N)$ (or $\calO(\log s)$ when $N\asymp s$, see
\eqref{PDE_Ns}). However, this method requires us to store all the vectors
$\boldsymbol{b}_{k, \ell}$. Using the sparsity of the stiffness matrices
which is of $\mathcal{O}(M)$, we require $\mathcal{O}(M N)$ storage.

\subsection{Partial differential equations with ``log-normal'' random
coefficients } \label{num_exp:3}
We consider the PDE \eqref{PDE_uniform} again but now we assume that the
random diffusion coefficient is given in parametric form by
\begin{equation*}
 \mathfrak{a}(\vec{x}, \boldsymbol{y})
= \exp\left( \psi_0(\vec{x}) + \sum_{j=1}^\infty y_j\, \psi_j(\vec{x}) \right),
 \quad y_j \sim \mbox{i.i.d. } \mathcal{N}(0,1).
\end{equation*}
This formula arises from the assumption that the logarithm of the random
coefficient $\mathfrak{a}(\vec{x},\cdot)$ is a Gaussian random field in
the domain $\dom$, which is parametrized in terms of principal components
of its covariance operator by a Karhunen-Lo\`eve expansion. We shall refer
to this case as the ``log-normal" case.

We may proceed as in the previous subsection, following
\eqref{PDE_lin_sys}--\eqref{eq_sol_PDE}. However, unlike the uniform case
where linearity can be exploited, the integral \eqref{a_approx} for the
log-normal case generally cannot be solved explicitly so that we need to
use a quadrature rule to approximate \eqref{a_approx}. Let $\vec{x}_{1,k,\ell},
\vec{x}_{2, k, \ell}, \ldots, \vec{x}_{I, k, \ell} \in \dom$ denote the set of quadrature points
and $w_{1, k, \ell}, w_{2, k, \ell}, \ldots, w_{I,k, \ell} \in\bbR$ denote the corresponding quadrature
weights which are used to approximate \eqref{a_approx},
\begin{equation}\label{def:b nkl}
  b_{n,k,\ell} \approx \widehat{b}_{n,k,\ell}
  = \sum_{i=1}^I w_{i, k, \ell } \, \mathfrak{a}(\vec{x}_{i, k, \ell}, \bsy_n) \nabla \phi_k(\vec{x}_{i, k, \ell}) \cdot \nabla \phi_\ell(\vec{x}_{i, k, \ell })
\end{equation}
for $0\le n < N$ and $1\le k,\ell\le M$. Let
$\widehat{B}(\boldsymbol{y}_n) = (\widehat{b}_{n,k,\ell})_{k, \ell}$. Thus
we need to compute
\begin{equation}\label{eq_as}
 \mathfrak{a}(\vec{x}_{i, k, \ell} , \boldsymbol{y}_n) = \exp (\theta_{i, k, \ell, n}),
 \quad
 \theta_{i, k, \ell, n}
 = \psi_0(\vec{x}_{i, k, \ell }) + \sum_{j=1}^s y_{n,j}\, \psi_j(\vec{x}_{i, k, \ell}),
\end{equation}
for all $1 \le i \le I$, $0 \le n < N$ and $1 \le k, \ell \le M$ such that $ \nabla \phi_k(\vec{x}_{i, k, \ell}) \cdot \nabla \phi_\ell(\vec{x}_{i, k, \ell} )$ is nonzero. The number of these nonzero inner products is $\calO(M)$ since for a fixed $k$, the number of $\ell$ such that the intersection of the supports of $\phi_k$ and $\phi_\ell$ is nonempty does not depend on $M$. The standard approach to obtain $\widehat{B}(\boldsymbol{y}_n)$ for all
$n=0,\ldots,N-1$ therefore requires $\calO(I\, M\, N\, s)$ operations.

We now describe the fast approach. Let
\begin{align*}
\Theta_{k, \ell} = & (\theta_{i, k, \ell, n})_{\substack{1 \le i \le I \\ 0 \le n < N}}, \\
\widehat{\Psi}_{i, k, \ell} = & \begin{pmatrix} \psi_0(\vec{x}_{i, k, \ell}) \\
\psi_0(\vec{x}_{i, k, \ell}) \\ \vdots \\ \psi_0(\vec{x}_{i, k, \ell}) \end{pmatrix},
\quad \widehat{\Psi}_{k, \ell} =( \widehat{\Psi}_{1, k, \ell}, \ldots, \widehat{\Psi}_{I, k, \ell}) \in \mathbb{R}^{N \times I}, \\
\Psi_{i, k, \ell} = & \begin{pmatrix} \psi_1(\vec{x}_{i, k, \ell}) \\ \psi_2(\vec{x}_{i, k, \ell} ) \\ \vdots \\ \psi_s(\vec{x}_{i, k, \ell}) \end{pmatrix},
\quad  \Psi_{k, \ell} = (\Psi_{1, k, \ell}, \ldots, \Psi_{I, k, \ell}) \in \mathbb{R}^{s \times I}.
\end{align*}
Then \eqref{eq_as} can be written in matrix form as
\begin{equation}\label{fastB}
 \Theta_{k, \ell} = \widehat{\Psi}_{k, \ell} +  Y \Psi_{k, \ell},
\end{equation}
where the multiplication $Y \Psi_{k, \ell}$ should be done as described in Section~\ref{sec_fast}. Hence \eqref{eq_as} can be computed using the fast QMC matrix method and
$\widehat{B}(\boldsymbol{y}_n)$ for all $0 \le n < N$ can be computed in
$\mathcal{O}(I\,M\, N \log N)$ operations. Again, the saving is that we replaced the factor $\mathcal{O}(s)$ by a factor $\mathcal{O}(\log N)$ in the computational cost.

\section{Numerical experiments} \label{sec:num}

In this section we carry out numerical experiments for the three
applications from the previous section. In all the numerical experiments
in the paper the times are averaged from $5$ independent runs using Matlab
R2013b on an Intel\textregistered Core$^{\tiny \mathrm{TM}}$ Xeon E5-2650v2
CPU {@} 2.6GHz.

\subsection*{Experiment 1: normally distributed points}

We are interested in comparing the computation times using the standard
approach of multiplying $\boldsymbol{y}_n A$ for $n = 0, 1,\ldots, N-1$,
and the fast QMC matrix approach described in
Subsection~\ref{sec_fast_lattice} with lattice point sets.

Table~\ref{table7} shows the computation times in seconds for various
values of $N$ and $s$. The value on top shows the standard approach,
whereas the value below shows the fast QMC matrix approach. For our
experiments we chose $\boldsymbol{\mu} = (0, 0, \ldots, 0)$, and $A$ a
random upper triangular matrix with positive diagonal entries (so that $A$
corresponds to the Cholesky factor of a random matrix $\Sigma$). The
computation times do not include the component-by-component construction
of the lattice generating vectors, the computation of $c_1,\ldots, c_s$
(since this information can be obtained from the fast
component-by-component construction), nor the computation of
$\Phi^{-1}(n/N)$ for $n = 0, 1, \ldots, N-1$ (since this computation is
the same for both methods).

The numerical experiments in Table~\ref{table7} show that there is an
advantage using the fast QMC matrix-vector product if the dimension is
large and the advantage grows as the dimension increases. This is in
agreement with the theory since the computational cost in the standard
approach is of order $\mathcal{O}(Ns^2)$ operations, whereas in the fast
QMC matrix-vector approach it is of order $\mathcal{O}(s\,N \log N)$
operations. Recall that the fast QMC matrix method incurs a storage cost
of $\mathcal{O}(Ns)$.

\begin{table}
\centering
\begin{tabular}{|c | c | c | c | c | c | c| }
\hline
Method & $N$ & $s=200$ & $s=400$ & $s=600$ & $s=800$ & $s=1000$   \\
\hline
std. & 16001    &  0.309 &   0.741 &  1.296 &    1.617 &    2.154 \\
fast &          &  0.164 &   0.301 &  0.450 &    0.589 &    0.741  \\ \hline
std. & 32003    &  0.589 &    1.468 & 2.435 & 3.063 & 4.238  \\
fast &          &  0.603 &    1.198 & 1.792 & 2.395 & 2.994  \\ \hline
std. & 64007    &  1.167 & 2.970 & 4.921 & 6.001 &  8.349   \\
fast &          &  1.804 & 3.853 & 5.551 & 7.582 &  9.827 \\ \hline
std. & 127997   &  2.579 &   5.889 &   9.490 &  11.891 &  16.818 \\
fast &          &  2.331 &   4.661 &   7.321 &  9.984  &  12.284  \\ \hline
std. & 256019   &  4.279  & 11.105 &  17.646 &  23.115 &  33.541   \\
fast &          &  5.401  & 10.933 &  16.174 &  24.147 &  26.898 \\ \hline
std. & 512009   & 8.885 &   23.368 &   31.942 &   48.059 &   66.378 \\
fast &          & 10.947 &   22.066 &   35.543 &  45.164  & 56.190 \\ \hline
\end{tabular} \vspace{0.3cm}
\caption{\small Times (in seconds) to generate normally distributed points with random covariance matrix.
The top row is the time required by using the standard approach, whereas the bottom row shows the time required using the fast QMC matrix-vector approach.}
\label{table7}
\end{table}

\subsection*{Experiment 2: the uniform case}

We consider the ODE
\begin{equation}\label{ode uniform}
\begin{aligned}
 -\frac{\mathrm{d}}{\mathrm{d} x} \left(a(x, \boldsymbol{y}) \frac{\mathrm{d}}{\mathrm{d} x} u(x, \boldsymbol{y}) \right)
 &= g(x) \quad\mbox{for } x \in (0,1) \mbox{ and } \boldsymbol{y} \in [-\tfrac{1}{2},\tfrac{1}{2}]^\mathbb{N}, \\
 u(x, \boldsymbol{y}) &= 0 \quad\mbox{for } x = 0, 1, \\
 a(x, \boldsymbol{y}) &= 2 + \sum_{j=1}^\infty y_j\, j^{-3/2} \sin(2\pi j x).
\end{aligned}
\end{equation}
Thus $\sum_{j\ge 1} \|\psi_j\|_{L^\infty(0,1)}^{2/3 + \varepsilon}<\infty$
for any $\varepsilon> 0$, and \eqref{PDE_Ns} implies that we should choose
$N \asymp s$. In our experiments we choose $M = N = s$.

To obtain an approximation of the solution we use finite elements. Let
$x_k = k/M$ for $k = 0,1, \ldots, M$ and for $k = 1, 2, \ldots, M-1$
define the hat function
\begin{equation}\label{def:phi 1d}
\phi_k(x) = \begin{cases} (x-x_{k-1}) M & \mbox{if } x_{k-1} \le x \le x_k, \\ (x_{k+1} - x) M & \mbox{if } x_k \le x \le x_{k + 1}, \\ 0 & \mbox{otherwise.}  \end{cases}
\end{equation}
Then
\begin{align*}
a_{0, k, \ell} = \begin{cases} 4 M & \mbox{if } k = \ell, \\ - 2 M & \mbox{if } |k-\ell| = 1, \\ 0 & \mbox{otherwise}, \end{cases}
\end{align*}
and for $j \ge 1$ we have
\begin{align*}
a_{j,k,\ell} = & \int_0^1 j^{-3/2} \sin(2\pi j x) \phi'_k(x) \phi'_\ell(x) \,\mathrm{d} x \\ = & \begin{cases} \frac{M^2}{\pi j^{5/2} } \sin \left( \frac{2\pi j}{M} \right) \sin \left( \frac{2\pi j k}{M} \right)  & \mbox{if } k = \ell, \\ - \frac{M^2}{\pi j^{5/2} }  \sin \left( \frac{\pi j}{M} \right) \sin \left( \frac{\pi j (2k-1)}{M} \right)  & \mbox{if } \ell = k-1, \\ - \frac{M^2}{\pi j^{5/2} }  \sin \left( \frac{\pi j}{M} \right) \sin \left( \frac{\pi j (2k + 1)}{M} \right)  & \mbox{if } \ell = k+1, \\ 0 & \mbox{otherwise}. \end{cases} 
\end{align*}
Thus the matrices $A_j$ and $B(\boldsymbol{y}_n)$ are tridiagonal. For
simplicity we choose $g$ such that
$(\widehat{g}_1,\widehat{g}_2,\ldots,\widehat{g}_M) = (1,1, \ldots, 1)$.

Table~\ref{table-new1} shows the computation times comparing the
standard approach with the fast QMC matrix method based on lattice point
sets as described in Section~\ref{sec_fast_lattice}. In this case the
mapping in $\bsy_n = \varphi(\bsx_n)$ is $\varphi(x) = x - 1/2$, since the
lattice points need to be translated from the usual unit cube $[0,1]^s$ to
$[-\tfrac{1}{2},\tfrac{1}{2}]^s$. Note that we do not apply any random
shifting as analyzed in \cite{KSS12}. Since the dimension $s$ is large,
the fast QMC matrix method is very effective in reducing the computation
times.
Note that in Table~\ref{table-new1} for the case $M=s=2N$ the times
for the standard method for $N=8009$ and $N=16001$ are in hours and are estimated from
extrapolating on previous values in the table. The experiments show there
is a clear advantage of fast QMC matrix-vector approach especially for
large values of $M,N$ and~$s$.

\begin{table}[t]
\centering
\begin{tabular}{| c | c | c | c | c | c | c | c | c| c |}
\hline
 \multicolumn{10}{|c|}{$M = s = 2N$} \\
 \hline
 $N$  &      67 &    127 &     257 &     509 &    1021 &    2053 &   4001 & 8009 & 16001  \\
\hline
std.  &   1     &  5     &     31  &     190 &   1346  &   10610 &  74550 & $\approx$144h
   & $\approx$1000h  \\
fast  &  0.035  &  0.042 &   0.114 &   0.462 &   1.562 &   5.591 & 19.678 & 87.246 &  342.615  \\
\hline
\hline
 \multicolumn{10}{|c|}{$M = s = \lceil\sqrt{N}\rceil$} \\
 \hline
 $N$  &      67 &    127 &     257 &     509 &    1021 &    2053 &   4001 & 8009 & 16001  \\
\hline
std.  &   0.066     &  0.164     &     0.474  &    1.272 &   3.570  &   10.813 &  30.127 & 89.42 & 273.873  \\
fast  &  0.012  &  0.015 &   0.028 &   0.059 &   0.126 &   0.265 & 0.516 & 1.113 &  2.443  \\
\hline
\hline
 \multicolumn{10}{|c|}{$s=N$ and $M=N^2$} \\
\hline
 $N$  &      67 &    127 &     257 &  509 & \multicolumn{5}{c|}{}   \\
\hline
std.  & 6    &  82  &  1699 & 27935 & \multicolumn{5}{c|}{} \\
fast  &0.243 & 1.385 & 11.268 &   107.042 & \multicolumn{5}{c|}{} \\
\hline
\end{tabular} \vspace{0.3cm}
\caption{\small Times (in seconds) to obtain the average value of the
finite element coefficients of the approximation \eqref{eq_sol_PDE}
to~\eqref{ode uniform}. Top: $M = s = 2N$. Middle: $M = s =
\lceil\sqrt{N}\rceil$. Bottom: $s=N$ and $M=N^2$.} \label{table-new1}
\end{table}

\subsection*{Experiment 3: the log-normal case}

In one space dimension, we consider the two-point boundary value problem for
the parametric, second order ODE
\begin{equation}\label{ode log}
\begin{aligned}
 -\frac{\mathrm{d}}{\mathrm{d} x} \left(a(x, \boldsymbol{y}) \frac{\mathrm{d}}{\mathrm{d} x} u(x, \boldsymbol{y}) \right)
 &= g(x) \quad\mbox{for } x \in (0,1) \mbox{ and } \boldsymbol{y} \in [-\tfrac{1}{2},\tfrac{1}{2}]^\mathbb{N}, \\
 u(x, \boldsymbol{y}) &= 0 \quad\mbox{for } x = 0, 1, \\
 a(x, \boldsymbol{y}) &= \exp\left(2 + \sum_{j=1}^\infty y_j\, j^{-3/2} \sin(2\pi j x)\right).
\end{aligned}
\end{equation}

We use $M$ finite elements to construct the approximate solutions as in \eqref{def:phi 1d}.
To compute \eqref{def:b nkl}, we use an equal weight quadrature with $M$ (so $I = M$) points. 

Table~\ref{table-new2} shows the computation time for the log-normal
case with different choices of number of finite elements $M$, the number
of QMC points $N$ and the truncated dimension $s$. As one would expect
from the theory, the most significant advantage of the fast QMC matrix
method occurs when $2^s$ is large compared to $N$, which is also reflected
in the numerical results.

\begin{table}[t]
\centering
\begin{tabular}{| c | c | c | c | c | c | c | c | c| c |}
\hline
 \multicolumn{10}{|c|}{$M = s = 2N$} \\
 \hline
 $N$  &      67 &    127 &     257 &     509 &    1021 &    2053 &   4001 & 8009 & 16001  \\
\hline
std.  &  0.028      &  0.051     &     0.140  &    0.436 &   1.734  & 15.173   & 84.381  & 614.636 & 4391.2   \\
fast  &  0.040  & 0.033  & 0.094   &  0.326  & 1.122   & 4.296   & 15.203  & 60.546 &  270.691  \\
\hline
\hline
 \multicolumn{10}{|c|}{$M = s = \lceil\sqrt{N}\rceil$} \\
 \hline
 $N$  &      67 &    127 &     257 &     509 &    1021 &    2053 &   4001 & 8009 & 16001  \\
\hline
std.  &  0.030  &  0.053 &   0.090 &   0.182 &   0.375 &   0.791 &   1.609 &   4.100 &   7.874 \\
fast & 0.132  &   0.036  &  0.052 &    0.106 &   0.228 &   0.480 &   0.940 &   2.670 &   4.597 \\
\hline
\hline
 \multicolumn{10}{|c|}{$s=N$ and $M=N^2$} \\
\hline
 $N$  &      67 &    127 &     257 &     509 & 1021 & \multicolumn{4}{c|}{}  \\
\hline
std.  & 0.162 &    0.945 &    9.935 &  84.790  & 891.175 & \multicolumn{4}{c|}{}    \\
fast  & 0.204 &   1.084  &  10.154  &  83.861  & 746.907 & \multicolumn{4}{c|}{}   \\
\hline
\end{tabular} \vspace{0.3cm}
\caption{\small Times (in seconds) to obtain the average value of the
finite element coefficients of the approximation \eqref{eq_sol_PDE} to
\eqref{ode log}. Top: $M = s = 2N$. Middle: $M = s =
\lceil\sqrt{N}\rceil$. Bottom: $s=N$ and $M=N^2$.} \label{table-new2}
\end{table}

\section*{Acknowledgment}
Josef Dick is the recipient of an Australian Research Council Queen
Elizabeth II Fellowship (project number DP1097023). Frances Y. Kuo is the
recipient of an Australian Research Council Future Fellowship (project
number FT130100655). Quoc T. Le Gia was supported partially by the ARC
Discovery Grant DP120101816. This research was supported under Australian Research Council's Discovery Projects funding scheme (project number DP150101770). Christoph Schwab acknowledges support through
ERC and SNF.

\bigskip

\noindent {\bf Addresses:} \\

Josef Dick, Frances Y. Kuo, Quoc T. Le Gia, School of Mathematics and Statistics, The University of New South Wales, Sydney, 2052 NSW, Australia. e-mail: {josef.dick, f.kuo, qlegia}(AT)unsw.edu.au \\

Christoph Schwab, Seminar for Applied Mathematics, ETH, 8092 Z\"urich, Switzerland. e-mail: christoph.schwab(AT)sam.math.ethz.ch


\begin{thebibliography}{99}

\bibitem{AP05} Y. Achdou and O. Pironneau, Computational methods for option pricing. Frontiers in Applied Mathematics, 30. Society for Industrial and Applied Mathematics (SIAM), Philadelphia, PA, 2005.

\bibitem{D14} J.~Dick,  Numerical integration of H\"older continuous, absolutely convergent Fourier, Fourier cosine, and Walsh series. J. Approx. Theory 183, 14--30 (2014).

\bibitem{DKS13}
  J.~Dick, F.~Y. Kuo, and I.~H. Sloan, High dimensional integration -- the quasi-Monte Carlo way. Acta Numer.\ 22, 133--288 (2013).

\bibitem{DNP13} J. Dick, D. Nuyens, and F. Pillichshammer, Lattice rules for non-periodic smooth integrands. Numer. Math. 126, 259--291 (2014).

\bibitem{DP07} J. Dick and F. Pillichshammer, Strong tractability of multivariate integration of arbitrary high order using digitally shifted polynomial lattice rules. J. Complexity 23, 436--453 (2007).

\bibitem{DP14} J. Dick and F. Pillichshammer, The weighted star-discrepancy of Korobov's $p$-sets. To appear in Proc. Amer. Math. Soc., 2015.

\bibitem{DP10} J. Dick and F. Pillichshammer, Digital nets and sequences. Discrepancy theory and quasi-Monte Carlo integration. Cambridge University Press, Cambridge, 2010.

\bibitem{DN97} C. R.~Dietrich and G. H.~Newsam, Fast and exact simulation
    of stationary Gaussian processes through circulant
    embedding of the covariance matrix. SIAM J. Sci. Comput. 18, 1088--1107 (1997).

\bibitem{FFTW05} M. Frigo and S.~G.~Johnson, The design and implementation of {FFTW3}. Proceedings of the IEEE 93, 216--231 (2005).

 
\bibitem{GKSW08} M. Giles, F. Y. Kuo, I. H. Sloan, and B. J. Waterhouse, Quasi-Monte Carlo for finance applications. ANZIAM Journal 50 (CTAC2008), C308 -- C323 (2008).

\bibitem{Gla} P. Glasserman, Monte Carlo methods in financial engineering. Applications of Mathematics (New York), 53. Stochastic Modelling and Applied Probability. Springer-Verlag, New York, 2004.

\bibitem{GKNSSS13} I.~G. Graham, F.~Y. Kuo, J.~A. Nichols, R.~Scheichl,
    Ch.~Schwab, and I.~H. Sloan, Quasi-Monte Carlo finite element methods for elliptic pdes with log-normal random coefficients. To appear in Numerische Math., 2015.

\bibitem{GKNSS11} I. G. Graham, F. Y. Kuo, D. Nuyens, R. Scheichl, and I. H. Sloan, Quasi-Monte Carlo methods for elliptic PDEs with random coefficients and applications. J. Comput. Phys. 230, 3668--3694 (2011).

\bibitem{HLD} W. H\"ormann, J. Leydold, and G. Derflinger, Automatic nonuniform random variate generation. Springer, Berlin, 2004.


\bibitem{HW81} L. K. Hua and Y. Wang, Applications of Number Theory to Numerical Analysis. Springer, Berlin, 1981.

\bibitem{KDSWW08} F.~Y.~Kuo, W.~T.~M.~Dunsmuir, I.~H.~Sloan, M.~P.~Wand, and R.~S.~Womersley, Quasi-Monte Carlo for highly structured generalised response models. Methodol.\ Comput.\ Appl.\ Probab.\ 10, 239--275 (2008).

\bibitem{KSS12} F. Y. Kuo, Ch. Schwab, and I. H. Sloan, Quasi-Monte Carlo finite element methods for a class of elliptic partial differential equations with random coefficients. SIAM J. Numer. Anal. 50, 3351--3374 (2012).

\bibitem{Leo12}  G. Leobacher, Fast orthogonal transforms and generation of Brownian paths. J. Complexity 28, 278--302 (2012).

\bibitem{NC06a} D. Nuyens and R. Cools, Fast algorithms for component-by-component construction of rank-1 lattice rules in shift-invariant reproducing kernel Hilbert spaces. Math. Comp. 75, 903--920 (2006).

\bibitem{NC06b} D. Nuyens and R. Cools, Fast component-by-component construction of rank-1 lattice rules with a non-prime number of points. J. Complexity 22, 4--28 (2006).

\bibitem{NC06c} D. Nuyens and R. Cools, Fast component-by-component construction, a reprise for different kernels. In: H. Niederreiter and D. Talay (eds.),  Monte Carlo and quasi-Monte Carlo methods 2004, 373--387, Springer, Berlin, 2006.

\bibitem{Sche07} K. Scheicher, Complexity and effective dimension of discrete L\'evy areas. J. Complexity 23, 152--168 (2007).

\bibitem{CSMCQMC12} Ch.~Schwab, QMC Galerkin discretization of parametric operator equations. In: J. Dick, F. Y. Kuo, G. W. Peters and I. H. Sloan (eds.), Monte Carlo and quasi-Monte Carlo methods 2012, 613--629, Springer, Berlin 2013.

\end{thebibliography}
\end{document}